\documentclass[12pt]{article}\usepackage{amsmath}

\usepackage{amssymb}

\begin{document}

\begin{titlepage}
\title{\bf Dynamical Systems on Bundles}
\author{ Mehmet Tekkoyun \footnote{E-mail address: tekkoyun@pau.edu.tr; Tel: +902582953616; Fax: +902582953593}\\
{\small Department of Mathematics, Pamukkale University,}\\
{\small 20070 Denizli, Turkey}}
\date{\today}
\maketitle
\begin{abstract}
In this study, Hamiltonian and Lagrangian theories, which are
mathematical models of mechanical systems, are structured on the
horizontal and the vertical distributions of tangent and cotangent
bundles. In the end, the geometrical and physical results related
to  Hamiltonian and Lagrangian dynamical systems are concluded.
\end{abstract}
\end{titlepage}

\section{Introduction}

The tangent bundle $TM$ and cotangent bundle $T^{*}M$ of manifold $M$ are in
good condition to be phase-spaces of velocity and momentum of a given
configuration space. Then modern differential geometry explains the
Lagrangians and Hamiltonian theories in classical mechanics. For example,
the tangent bundle $TM$\ carries some natural object fields, as: Liouville
vector field $V$, tangent structure $J$, almost product structure $P$, the
vertical distribution $V$, the horizontal distribution $H$, semispray $X$.
Therefore these structures have an important role in physical fields. The
following symbolical equation expresses the dynamical equations for both
theories:
\begin{equation}
i_{X}\Phi =\digamma  \label{1.1}
\end{equation}
If one studies the Hamiltonian theory then equation (\ref{1.1}) is the
intrinsical form of Hamiltonian equations, where $\Phi =\phi _{\mathbf{H}%
}=-d\lambda $ , $\digamma =d\mathbf{H}$ and $\lambda $ is Liouville form
canonically constructed on the cotangent bundle $T^{*}M$ of the
configuration manifold $M$, and $\mathbf{H}$ is a $C^{\infty }-$ function on
$T^{*}M$ such that $\mathbf{H}:T^{*}M\rightarrow \mathbf{R.}$ If one studies
the Lagrangian theory then equation (\ref{1.1}) is the intrinsical form of
Lagrangian equations, where $\Phi =\Phi _{L}=-dd_{P}L$ such that $%
L:TM\rightarrow \mathbf{R}$ Lagrangian function, $\digamma =dE_{L}$, $%
E_{L}=V(L)-L\mathbf{,}$ $V$ the Liouville vertical vector field on $TM$ , $%
E_{L}$ the energy associated to the function $L$ and $X$ is Liouville
vertical vector field canonically constructed on the tangent bundle $TM$ of
the configuration manifold $M.$ Mathematical expressions of mechanical
systems are always determined by the Hamiltonian and Lagrangian systems.
These expressions, in particular geometric expressions in mechanics and
dynamics, are given in some studies \cite{mcrampin, nutku,deleon,deleon1}.
Klein introduced that the geometric study of the Lagrangian theory admits an
alternative approach without the use of the regular condition on $L$ \cite%
{klein}. Para-complex analogues of the Lagrangians and Hamiltonians were
obtained in the framework of K\"{a}hlerian manifold and the geometric
conclusions on a para-complex dynamical systems were obtained \cite%
{tekkoyun2005}. As well-known from before works, Lagrangian distribution on
symplectic manifolds are used in geometric quantization and a connection on
a symplectic manifold is an important structure to obtain a deformation
quantization \cite{etayo}. In the before studies; although
real/(para)complex geometry and real/(para)complex mechanical-dynamical
systems were analyzed successfully, they have not dealt with dynamical
systems on horizontal distribution $HTM$\ and vertical distribution $VTM$ of
tangent bundle $TM$ of manifold $M$. In this Letter, therefore,
Euler-Lagrange and Hamiltonian equations related to dynamical systems on the
distributions used in obtaining geometric quantization have been given.

\section{Preliminaries}

In this paper all geometrical object fields and all mappings are considered
of the class $C^{\infty },$\ expressed by the words ''differentiate'' or
''smooth''. \ The indices $i$,$j$,...run over set $\{1,..,n\}$\ and Einstein
convention of summarizing is adopted over all this paper. $\mathbf{R}$, $%
\mathcal{F}(TM)$, $\chi (TM),\chi (T^{*}M)$ denote the set of real numbers,
the set of real functions on $TM$, the set of vector fields on $TM$ and the
set of 1-forms on $T^{*}M$.

\subsection{Manifold, Bundle and Distributions}

In this subsection, some definitions were derived and taken from \cite{radu}%
. Let $TM$ be tangent bundle of a real $n$-dimensional differentiable
manifold $M$. Then it will be denoted a point of $M$ by $x$ and its local
coordinate system by $(U,\varphi )$ such that $\varphi (x)=(x^{i}).$\ Such
that the projection $\pi :TM\rightarrow M$, $\pi (u)=x,$\ a point $u\in TM$
will be denoted by $(x,y)$, its local coordinates being $(x^{i},y^{i})$.
There are the natural basis $(\frac{\partial }{\partial x^{i}},\frac{%
\partial }{\partial y^{i}})$ and dual basis $(dx^{i},dy^{i})$\ of the
tangent space $T_{u}TM$\ and the cotangent space $T_{u}^{*}(TM)$ at the
point $u\in TM$\textbf{, }respectively$.$\ Consider the $\mathcal{F}(TM)-$
and $\mathcal{F}(T^{*}M)-$ linear mappings (also named to be almost tangent
structures) $J:\chi (TM)\rightarrow \chi (TM)$\thinspace and $J^{*}:\chi
(TM)\rightarrow \chi (TM)$ given by
\[
\begin{array}{l}
J(\frac{\partial }{\partial x^{i}})=\frac{\partial }{\partial y^{i}},\,J(%
\frac{\partial }{\partial y^{i}})=0,\,%
\end{array}
\]
and
\[
\begin{array}{l}
J^{*}(dx^{i})=dy^{i},J^{*}(dy^{i})=0.%
\end{array}
\]
The tangent space $V_{u}$ to the fibre $\pi ^{-1}(x)$\ in the point $u\in TM$%
\ is locally spanned by $\{\frac{\partial }{\partial y^{1}},..,\frac{%
\partial }{\partial y^{n}}\}$. Therefore, the mapping $V$: $u\in
TM\rightarrow V_{u}\subset T_{u}TM$\ provides a regular distribution
generated by the adapted basis $\{\frac{\partial }{\partial y^{i}}\}.$\
Consequently, $V$\ is an integrable distribution on $TM$. $V$\ is called the
vertical distribution on $TM$.\ Let $N$ be a nonlinear connection on $TM$. $%
N $ is characterized by $v$, $h$ vertical and horizontal projectors. We
consider the vertical projector $v:\chi (TM)\rightarrow \chi (TM)$ defined
by $v(X)=X,\,\forall \,X\in \chi (VTM);$ $v(X)=0,\forall \,X\in \chi (HTM).$
Similarly, the mapping $H$: $u\in TM\rightarrow H_{u}\subset T_{u}TM$\
provides a regular distribution determined by the adapted basis $\{\frac{%
\delta }{\delta x^{i}}\}.$\ Consequently, $H$\ is an integrable distribution
on $TM$. $H$\ is called the horizontal distribution on $TM$.\ There is a $%
\mathcal{F}(TM)-$linear mapping $h:\chi (TM)\rightarrow \chi (TM),$ for
which $h^{2}=h,$ $Ker$ $h=\chi (VTM).$ Any vector field $X\in \chi (TM)$ can
be uniquely written as follows $X=hX+vX=X^{H}+X^{V}$. Therefore $X^{H}$ and $%
X^{V}$ are horizontal and vertical components of vector field $X$.
Therefore, any vector field $X$ can be uniquely written in the form
\[
\begin{array}{l}
X=X^{H}+X^{V}%
\end{array}
\]
such that
\[
\begin{array}{l}
X^{H}=X^{i}(\frac{\partial }{\partial x^{i}}-N_{j}^{i}(x,y)\frac{\partial }{%
\partial y^{j}}),\,\,\,X^{V}=X^{i}N_{j}^{i}(x,y)\frac{\partial }{\partial
y^{j}}%
\end{array}
\]
where $N_{j}^{i}$ are local coefficients of nonlinear connection $N$ on $TM$
$P$ is an almost product structure on $TM.$ $(\frac{\delta }{\delta x^{i}},%
\frac{\partial }{\partial y^{i}})$ is a local basis adapted to the
horizontal distribution $HTM$ and the vertical distribution $VTM$. Then $%
(dx^{i},\delta y^{i})$ is dual basis of $(\frac{\delta }{\delta x^{i}},\frac{%
\partial }{\partial y^{i}})$ basis. We have
\[
\begin{array}{l}
P(X)=X,\forall \,X\in \chi (HTM);\,\,\,P(X)=-X,\forall \,X\in \chi (VTM) \\
P^{*}(\omega )=\omega ,\forall \,\omega \in \chi
(HT^{*}M);\,\,\,P^{*}(\omega )=-\omega ,\forall \,\omega \in \chi (VT^{*}M),%
\end{array}
\]
where $P^{*}$ is the dual structure of $P$ For the $(\frac{\delta }{\delta
x^{i}},\frac{\partial }{\partial y^{i}})$ basis and $(dx^{i},\delta y^{i})$
dual basis we have
\[
\begin{array}{l}
\frac{\delta }{\delta x^{i}}=\frac{\partial }{\partial x^{i}}-N_{j}^{i}(x,y)%
\frac{\partial }{\partial y^{j}}.%
\end{array}
\]
and
\[
\begin{array}{l}
\delta y^{i}=dy^{i}+N_{j}^{i}(x,y)dx^{j}.%
\end{array}
\]
For the operators $h,v,P,P^{*},J,J^{*}$ we get
\[
\begin{array}{l}
h+v=I;\,\,P=2h-I;\,\,P=h-v;\,\,P=I-2v, \\
JP=J;PJ=-J;J^{*}P^{*}=J^{*};P^{*}J^{*}=-J^{*}, \\
h(\frac{\delta }{\delta x^{i}})=\frac{\delta }{\delta x^{i}};\,h(\frac{%
\partial }{\partial y^{i}})=0;\,v(\frac{\delta }{\delta x^{i}})=0;v(\frac{%
\partial }{\partial y^{i}})=\frac{\partial }{\partial y^{i}},\, \\
P(\frac{\delta }{\delta x^{i}})=\frac{\delta }{\delta x^{i}};\,\,P(\frac{%
\partial }{\partial y^{i}})=-\frac{\partial }{\partial y^{i}}, \\
P^{*}(dx^{i})=dx^{i};\,P^{*}(\delta y^{i})=-\delta y^{i}.%
\end{array}
\]

\section{Lagrangian Dynamical Systems}

In this section, Euler-Lagrange equations for classical mechanics are
structured by means of almost product structure $P$\ under the consideration
of the basis $\{\frac{\delta }{\delta x^{i}},\,\frac{\partial }{\partial
y^{i}}\}$ on distributions $HTM$ and $VTM$ of tangent bundle $TM$ of
manifold $M.$ Let $(x^{i},y^{i})$ be its local coordinates. Let semispray be
the vector field $X$ given by
\begin{equation}
\begin{array}{l}
X=X^{i}\frac{\delta }{\delta x^{i}}+\overset{.}{X}^{i}\frac{\partial }{%
\partial y^{i}},\overset{.}{\,\,X}^{i}=X^{i}N_{j}^{i}%
\end{array}
\label{3.1}
\end{equation}
where the dot indicates the derivative with respect to time $t$. The vector
field denoted by $V=P(X)$ and expressed by
\begin{equation}
\begin{array}{l}
V=X^{i}\frac{\delta }{\delta x^{i}}-\overset{.}{X}^{i}\frac{\partial }{%
\partial y^{i}}%
\end{array}
\label{3.2}
\end{equation}
is called \textit{Liouville vector field} on the bundle $TM$. The maps given
by $\mathbf{T,P}:TM\rightarrow \mathbf{R}$ such that $\mathbf{T}=\frac{1}{2}%
m_{i}(x^{i})^{2},\mathbf{P}=m_{i}gh$ are called \textit{the kinetic energy}
and \textit{the potential energy of the mechanical system,} respectively.%
\textit{\ }Here\textit{\ }$m_{i},g$ and $h$ stand for mass of a mechanical
system having $m$ particles, the gravity acceleration and distance to the
origin of a mechanical system on the tangent bundle $TM$, respectively. Then
$L:TM\rightarrow \mathbf{R}$ is a map that satisfies the conditions; i) $L=%
\mathbf{T-P}$ is a \textit{Lagrangian function, ii)} the function given by $%
E_{L}=V(L)-L$ is \textit{a Lagrangian energy}. The operator $i_{P}$ induced
by $P$ and shown by
\begin{equation}
\begin{array}{l}
i_{P}\omega (X_{1},X_{2},...,X_{r})=\sum_{i=1}^{r}\omega
(X_{1},...,P(X_{i}),...,X_{r})%
\end{array}
\label{3.3}
\end{equation}
is said to be \textit{vertical derivation, }where $\omega \in \wedge
^{r}{}TM,$ $X_{i}\in \chi (TM).$ The \textit{vertical differentiation} $%
d_{P} $ is defined by
\begin{equation}
\begin{array}{l}
d_{P}=[i_{P},d]=i_{P}d-di_{P}%
\end{array}
\label{3.4}
\end{equation}
where $d$ is the usual exterior derivation. For an almost product structure $%
P$, the closed fundamental form is the closed 2-form given by $\Phi
_{L}=-dd_{P}L$ such that
\begin{equation}
\begin{array}{l}
d_{P}:\mathcal{F}(TM)\rightarrow {}T^{*}M%
\end{array}
\label{3.5}
\end{equation}
Then we have
\begin{equation}
\begin{array}{l}
\Phi _{L}=-(\frac{\delta }{\delta x^{j}}dx^{j}+\frac{\partial }{\partial
y^{j}}\delta y^{j})(\frac{\delta L}{\delta x^{i}}dx^{i}-\frac{\partial L}{%
\partial y^{i}}\delta y^{i}) \\
\,\,\,\,\,\,\,\,\,\,=\frac{\delta ^{2}L}{\delta x^{j}\delta x^{i}}%
dx^{j}\wedge dx^{i}-\frac{\delta (\partial L)}{\delta x^{j}\partial y^{i}}%
dx^{j}\wedge \delta y^{i}-\frac{\partial (\delta L)}{\partial y^{j}\delta
x^{i}}\delta y^{j}\wedge dx^{i}+\frac{\partial ^{2}L}{\partial y^{j}\partial
y^{i}}\delta y^{j}\wedge \delta y^{i}.%
\end{array}
\label{3.6}
\end{equation}
Let $X$ be the second order differential equation (semispray) determined by
\textbf{Eq. }(\ref{1.1}) and
\begin{equation}
\begin{array}{l}
i_{X}\Phi _{L}=\Phi _{L}(X)=-X^{i}\frac{\delta ^{2}L}{\delta x^{j}\delta
x^{i}}\delta _{i}^{j}dx^{i}+X^{i}\frac{\delta ^{2}L}{\delta x^{j}\delta x^{i}%
}dx^{j}+X^{i}\frac{\delta (\partial L)}{\delta x^{j}\partial y^{i}}\delta
_{i}^{j}\delta y^{i}-\overset{.}{X}^{i}\frac{\delta (\partial L)}{\delta
x^{j}\partial y^{i}}dx^{j} \\
-\overset{.}{X}^{i}\frac{\partial (\delta L)}{\partial y^{j}\delta x^{i}}%
\delta _{i}^{j}dx^{i}+X^{i}\frac{\partial (\delta L)}{\partial y^{j}\delta
x^{i}}\delta y^{j}+\overset{.}{X}^{i}\frac{\partial ^{2}L}{\partial
y^{j}\partial y^{i}}\delta _{i}^{j}\delta y^{i}-\overset{.}{X}^{i}\frac{%
\partial ^{2}L}{\partial y^{j}\partial y^{i}}\delta y^{j}.%
\end{array}
\label{3.7}
\end{equation}
Since the closed 2-form $\Phi _{L}$ on $TM$ is in the symplectic structure,
it is found
\begin{equation}
\begin{array}{l}
E_{L}=V(L)-L=X^{i}\frac{\delta L}{\delta x^{i}}-\overset{.}{X}^{i}\frac{%
\partial L}{\partial y^{i}}-L%
\end{array}
\label{3.8}
\end{equation}
and hence
\begin{equation}
\begin{array}{l}
dE_{L}=X^{i}\frac{\delta ^{2}L}{\delta x^{j}\delta x^{i}}dx^{j}-\overset{.}{X%
}^{i}\frac{\delta (\partial L)}{\delta x^{j}\partial y^{i}}dx^{j}-\frac{%
\delta L}{\delta x^{j}}dx^{j}+X^{i}\frac{\partial (\delta L)}{\partial
y^{j}\delta x^{i}}\delta y^{j}-\overset{.}{X}^{i}\frac{\partial ^{2}L}{%
\partial y^{j}\partial y^{i}}\delta y^{j}-\frac{\partial L}{\partial y^{j}}%
\delta y^{j}%
\end{array}
\label{3.9}
\end{equation}
With the use of \textbf{Eq.} (\ref{1.1}), considering (\ref{3.7}) and (\ref%
{3.9}) we get
\begin{equation}
\begin{array}{l}
-X^{i}\frac{\delta ^{2}L}{\delta x^{j}\delta x^{i}}\delta
_{i}^{j}dx^{i}+X^{i}\frac{\delta ^{2}L}{\delta x^{j}\delta x^{i}}dx^{j}+X^{i}%
\frac{\delta (\partial L)}{\delta x^{j}\partial y^{i}}\delta _{i}^{j}\delta
y^{i}-\overset{.}{X}^{i}\frac{\delta (\partial L)}{\delta x^{j}\partial y^{i}%
}dx^{j} \\
-\overset{.}{X}^{i}\frac{\partial (\delta L)}{\partial y^{j}\delta x^{i}}%
\delta _{i}^{j}dx^{i}+X^{i}\frac{\partial (\delta L)}{\partial y^{j}\delta
x^{i}}\delta y^{j}+\overset{.}{X}^{i}\frac{\partial ^{2}L}{\partial
y^{j}\partial y^{i}}\delta _{i}^{j}\delta y^{i}-\overset{.}{X}^{i}\frac{%
\partial ^{2}L}{\partial y^{j}\partial y^{i}}\delta y^{j} \\
=X^{i}\frac{\delta ^{2}L}{\delta x^{j}\delta x^{i}}dx^{j}-\overset{.}{X}^{i}%
\frac{\delta (\partial L)}{\delta x^{j}\partial y^{i}}dx^{j}-\frac{\delta L}{%
\delta x^{j}}dx^{j}+X^{i}\frac{\partial (\delta L)}{\partial y^{j}\delta
x^{i}}\delta y^{j}-\overset{.}{X}^{i}\frac{\partial ^{2}L}{\partial
y^{j}\partial y^{i}}\delta y^{j}-\frac{\partial L}{\partial y^{j}}\delta
y^{j}%
\end{array}
\label{3.10}
\end{equation}
or
\begin{equation}
\begin{array}{l}
-X^{i}\frac{\delta ^{2}L}{\delta x^{j}\delta x^{i}}dx^{j}-\overset{.}{X}^{i}%
\frac{\partial (\delta L)}{\partial y^{j}\delta x^{i}}dx^{j}+\frac{\delta L}{%
\delta x^{j}}dx^{j}+X^{i}\frac{\delta (\partial L)}{\delta x^{j}\partial
y^{i}}\delta y^{j}+\overset{.}{X}^{i}\frac{\partial ^{2}L}{\partial
y^{j}\partial y^{i}}\delta y^{j}+\frac{\partial L}{\partial y^{j}}\delta
y^{j}=0,%
\end{array}
\label{3.11}
\end{equation}
If a curve denoted by $\alpha :\mathbf{R}\rightarrow TM$ is considered to be
an integral curve of $X,$i.e$.$ $X(\alpha (t))=\frac{d\alpha (t)}{dt}$ then
\[
\begin{array}{l}
-\frac{d}{dt}(\frac{\delta L}{\delta x^{i}})+\frac{\delta L}{\delta x^{i}}+%
\frac{d}{dt}(\frac{\partial L}{\partial y^{i}})+\frac{\partial L}{\partial
y^{i}}=0.%
\end{array}
\]
or
\begin{equation}
\begin{array}{l}
\frac{d}{dt}(\frac{\delta L}{\delta x^{i}})-\frac{\partial L}{\partial y^{i}}%
=0,\,\,\frac{d}{dt}(\frac{\partial L}{\partial y^{i}})+\frac{\delta L}{%
\delta x^{i}}=0%
\end{array}
\label{3.12}
\end{equation}
Thus the equations\textbf{\ }given by\textbf{\ }(\ref{3.12}) are seen to be
a \textit{Euler-Lagrange equations} on $HTM$ horizontal and $VTM$ vertical
distributions$,$ and then the triple $(TM,\Phi _{L},X)$ is seen to be a
\textit{mechanical system }with taking into account the basis $\{\frac{%
\delta }{\delta x^{i}},\frac{\partial }{\partial y^{i}}\}$ on the
distributions $HTM$ and $VTM$.

\section{Hamiltonian Dynamical Systems}

In this section, Hamiltonian equations for classical mechanics are
structured on the distributions ${}HT^{*}M$ and $VT^{*}M$ of $T^{*}M.$
Suppose that an almost product structure, a Liouville form and a 1-form on $%
T^{*}M$ are shown by $P^{*}$, $\lambda $ and $\omega $, respectively$.$ Then
we hold
\begin{equation}
\begin{array}{l}
\omega =\frac{1}{2}(y^{i}dx^{i}+x^{i}\delta y^{i})%
\end{array}
\label{4.1}
\end{equation}
and
\begin{equation}
\begin{array}{l}
\lambda =P^{*}(\omega )=\frac{1}{2}(y^{i}dx^{i}-x^{i}\delta y^{i}).%
\end{array}%
\end{equation}
It is concluded that if $\phi $ is a closed 2- form on $T^{*}M,$ then $\phi
_{\mathbf{H}}$ is also a symplectic structure on ${}{}T^{*}M$. If
Hamiltonian vector field $X_{\mathbf{H}}$ associated with Hamiltonian energy
$\mathbf{H}$ is given by
\begin{equation}
\begin{array}{l}
X_{\mathbf{H}}=X^{i}\frac{\delta }{\delta x^{i}}+Y^{i}\frac{\partial }{%
\partial y^{i}},%
\end{array}
\label{4.3}
\end{equation}
then
\begin{equation}
\begin{array}{l}
\phi _{\mathbf{H}}=-d\lambda =-\delta y^{i}\wedge dx^{i}%
\end{array}%
\end{equation}
and
\begin{equation}
\begin{array}{l}
i_{X_{\mathbf{H}}}\phi =-Y^{i}dx^{i}+X^{i}\delta y^{i}.%
\end{array}
\label{4.5}
\end{equation}
Moreover, the differential of Hamiltonian energy is written as follows:
\begin{equation}
\begin{array}{l}
d\mathbf{H}=\frac{\delta \mathbf{H}}{\delta x^{i}}dx^{i}+\frac{\partial
\mathbf{H}}{\partial y^{i}}\delta y^{i}.%
\end{array}
\label{4.6}
\end{equation}
By means of \textbf{Eq.}(\ref{1.1}), using \textbf{Eq. }(\ref{4.5}) and
\textbf{Eq. }(\ref{4.6}), the Hamiltonian vector field is calculated to be
\begin{equation}
\begin{array}{l}
X_{\mathbf{H}}=\frac{\partial \mathbf{H}}{\partial y^{i}}\frac{\delta }{%
\delta x^{i}}-\frac{\delta \mathbf{H}}{\delta x^{i}}\frac{\partial }{%
\partial y^{i}}.%
\end{array}
\label{4.7}
\end{equation}
Suppose that a curve
\begin{equation}
\begin{array}{l}
\alpha {:\,}I\subset \mathbf{R}\rightarrow T^{*}M%
\end{array}
\label{4.8}
\end{equation}
be an integral curve of the Hamiltonian vector field $X_{\mathbf{H}}$, i.e.,
\begin{equation}
\begin{array}{l}
X_{\mathbf{H}}(\alpha (t))=\frac{d\alpha (t)}{dt},\,\,t\in I.%
\end{array}
\label{4.9}
\end{equation}
In the local coordinates, it is concluded that
\begin{equation}
\begin{array}{l}
\alpha (t)=(x^{i}(t),y^{i}(t))%
\end{array}
\label{4.10}
\end{equation}
and
\begin{equation}
\begin{array}{l}
\frac{d\alpha (t)}{dt}=\frac{dx^{i}}{dt}\frac{\delta }{\delta x^{i}}+\frac{%
dy^{i}}{dt}\frac{\partial }{\partial y^{i}}.%
\end{array}
\label{4.11}
\end{equation}
Taking into consideration \textbf{Eqs. }(\ref{4.8}), (\ref{4.6}),\textbf{\ }(%
\ref{4.10}), the result equations can be found
\begin{equation}
\begin{array}{l}
\frac{dx^{i}}{dt}=\frac{\partial \mathbf{H}}{\partial y^{i}},\frac{dy^{i}}{dt%
}=-\frac{\delta \mathbf{H}}{\delta x^{i}}%
\end{array}
\label{4.12}
\end{equation}
Thus, the equations\textbf{\ }(\ref{4.12}) are seen to be \textit{%
Hamiltonian equations} on the horizontal distribution ${}{}HT^{*}M$ and
vertical distribution $VT^{*}M.,$ and then the triple $(T^{*}M,\phi _{%
\mathbf{H}},X_{\mathbf{H}})$ is seen to be a \textit{Hamiltonian mechanical
system }with the use of basis $\{\frac{\delta }{\delta x^{i}},\frac{\partial
}{\partial y^{i}}\}$ on the distributions ${}HT^{*}M$ and $VT^{*}M$.

\section{Conclusions}

This paper has obtained to exist physical proof of the both mathematical
equality given by $TM=HTM\oplus VTM$ and its dual equality$.$ Lagrangian and
Hamiltonian dynamics have intrinsically been described with taking into
account the basis $\{\frac{\delta }{\delta x^{i}},\frac{\partial }{\partial
y^{i}}\}$ and dual basis $(dx^{i},\delta y^{i})$\ on distributions of
tangent and cotangent bundles $TM$ and $T^{*}M$ of manifold $M$..

\section{\textbf{Discussions}}

As known, geometry of Lagrangians and Hamiltonians give a model for
Relativity, Gauge Theory and Electromagnetism in a very natural blending of
the geometrical structures of the space with the characteristics properties
of these physical fields. Therefore we consider that the equations (\ref%
{3.12}) and (\ref{4.12}) especially can be used in fields determined the
above of physical.

\end{document}